\documentclass[10pt,a4paper,twoside]{amsart}

\usepackage{amsfonts,amssymb,euscript,eufrak}
\usepackage{amscd,amsmath}
\usepackage{cite}
\usepackage[all]{xy}

\newcommand{\bbA}{{\mathbb A}}
\newcommand{\bbB}{{\mathbb B}}
\newcommand{\bbC}{{\mathbb C}}

\newcommand{\bbQ}{{\mathbb Q}}

\newcommand{\bbZ}{{\mathbb Z}}

\newcommand{\Q}{\mathbb{Q}}
\newcommand{\C}{\mathbb{C}}

\newcommand{\Z}{\mathbb{Z}}
\newcommand{\N}{\mathbb{N}}

\newcommand{\ol}{o}

\newcommand{\gor}{gr^\bullet\,}

\newcommand{\G}{\mathbb{G}}

\newcommand{\Os}{\mathcal{O}}
\newcommand{\Fc}{\mathcal{F}}
\newcommand{\Gm}{\mathbb{G}_m}
\newcommand{\hol}{\hat{\ol}}\newcommand{\holn}{\hat{\ol}_n}\newcommand{\holnk}{\hat{\ol}_{K,n}}
\newcommand{\limi}{\underleftarrow{\lim}_n}

\newcommand{\B}{{\bf B}}

\newcommand{\hato}{\,\hat{\otimes}}

\newcommand{\Ft}{F_{t'}}

\newcommand{\ZK}{\hat{Z}_K}
\newcommand{\car}{\stackrel{\cong}{\longrightarrow}}

\newcommand{\Bc}{\mathring{B}}

\newcommand{\Kc}{\mathring{K}}

\newcommand{\Bcm}{\mathring{B}^{(m)}}

\newenvironment{pr}{\it Proof:\rm}{\hfill $\Box$\newline}
\newtheorem{theo}{Theorem}[section]
\newtheorem{prop}[theo]{Proposition}

\newtheorem{lem}[theo]{Lemma}

\begin{document}
\title[Picard groups]{Picard groups in $p$-adic Fourier theory}

\author[Tobias Schmidt]{Tobias Schmidt}
\address{Mathematisches Institut\\ Westf\"alische Wilhelms-Universit\"at
M\"unster\\ Einsteinstr. 62\\ D-48149 M\"unster, Germany}
\email{toschmid@math.uni-muenster.de}

\maketitle
\begin{abstract}
Let $L\neq\mathbb{Q}_p$ be a proper finite field extension of
$\bbQ_p$ and $\ol\subset L$ its ring of integers viewed as an
abelian locally $L$-analytic group. Let $\hol$ be the rigid
$L$-analytic group parametrizing the locally analytic characters
of $\ol$ constructed by Schneider-Teitelbaum. Let $K/L$ be a
finite extension field. We show that the base change $\hol_K$ has
a Picard group $Pic(\hol_K)$ which is profinite and that the unit
section in $\hol_K$ provides a divisor class of infinite order. In
particular, the abelian group $Pic(\hol_K)$ is not finitely
generated and is not a torsion group. On the way we show that
$\hol_K$ is a nontrivial \'etale covering of the affine line over
$K$ realized via the logarithm map of a Lubin-Tate formal group.
We finally prove that rank and determinant mappings induce an
isomorphism between $K_0(\hol_K) $ and $\bbZ \oplus Pic(\hol_K)$.
\end{abstract}
\section{Introduction}
\footnote[0]{2000 Mathematics Subject Classification: 22E50
(primary), 14G22 (secondary).} 

\bigskip

Let $L$ be a complete nonarchimedean field and $\B$ a rigid
$L$-analytic open polydisc. A twisted form of $\B$ is a rigid
analytic space $X$ over $L$ together with a complete
nonarchimedean extension $L\subseteq L'$ and an isomorphism
$X_{L'}\car \B_{L'}$. The study of the forms of $\B$ with respect
to a given extension $L\subseteq L'$ is an important yet difficult
problem which is still in its infancy. As a first result A. Ducros
recently has shown \cite{Ducros} that, if $L\subseteq L'$ is
finite and tamely ramified, any form $X$ is already trivial, in
the sense that $X$ is itself isomorphic to an open polydisc. He
also showed that there are finite wildly ramified extensions that
support plenty of nontrivial forms, even in dimension one.

\vskip8pt

A coarser but more accessible problem concerns the description of
basic invariants of such forms such as the Picard group or the
Grothendieck group. The aim of this note is to make a first step
in this direction and study the Picard group of an interesting
form that comes from $p$-adic representation theory. To give more
details, let $L\neq\mathbb{Q}_p$ be from now on a proper finite
extension of $\bbQ_p$ and let $\B$ be the open unit disc of
dimension one. In \cite{ST2} Schneider-Teitelbaum generalize the
classical $p$-adic Fourier theory of Y. Amice \cite{Amice2} for
the group $\bbZ_p$ to the additive group of integers $o$ in $L$.
The key step is the construction of a certain nontrivial form
$\hol$ of $\B$ ($\hol$ is even a group object) with respect to the
transcendental extension $L\subset \bbC_p$ having the surprising
property of admitting no trivializiation over any discretely
valued complete subfield of $\bbC_p$. In particular, this forces
its Picard group to be nontrivial. In \cite{Teitelbaum} J.
Teitelbaum suggested to study this Picard group, partly for
representation-theoretic reasons which we will describe below.

\vskip8pt

Slightly more general, we will study the series of Picard groups
$Pic(\hol_K)$ where $\hol_K$ denotes the base change of $\hol$ to
finite extensions $K$ of $L$. Although we are not able to
determine the structure of $Pic(\hol_K)$ explicitly, our main
result shows that it is of enormous size. More precisely, we prove
that $Pic(\hol_K)$ is a profinite group and that the unit section
of $\hol_K$ supports a divisor class of infinite order. In
particular, the abelian group $Pic(\hol_K)$ is not finitely
generated and is not a torsion group. On the way we show that
$\hol_K$ is a nontrivial \'etale covering of the affine line over
$K$ (in the sense of \cite{deJongET}). We finally show that rank
and determinant mappings induce an isomorphism $K_0(\hol_K)
\car\bbZ \oplus Pic(\hol_K)$.

\vskip8pt

As already indicated the variety $\hol_K$ has a
representation-theoretic interpretation. Indeed the additive group
$\ol$ is among the first examples of a compact abelian locally
$L$-analytic group. The generalized Amice-Fourier isomorphism of
Schneider-Teitelbaum \cite{ST2} induces an isomorphism of the ring
of holomorphic functions on $\hol_K$ with the $K$-valued locally
analytic distribution algebra of $o$. Since $\hol_K$ is a
quasi-Stein space in the sense of R. Kiehl \cite{Kiehl} the group
$Pic(\hol_K)$ controls the ideal structure of this distribution
algebra and therefore the locally analytic representation theory
of $o$. For example, the rational points of $\hol_K$ are in
bijective correspondence with the locally analytic characters
$o\rightarrow K^\times$.


\vskip8pt

In the following we briefly outline the article. As with most
nontrivial Picard groups there is by no means a straightforward
way to compute the structure of $Pic(\hol_K)$. Our strategy is to
first determine the local Picard groups corresponding to a
suitable open affinoid covering of the quasi-Stein space $\hol_K$
and then 'glue' these informations. We begin by recalling some
results on the divisor theory for Dedekind domains. The point is
that the affinoid algebra of a twisted form of a nonarchimedean
closed disc with respect to a finite extension is a Dedekind
domain. In sect. \ref{section-forms} we use descent theory to
prove a finiteness result for the Picard group of a certain class
of such forms. In sect. \ref{section-fourier} we turn to the
variety $\hol_K$. We prove that it admits an admissible affinoid
covering $\hol_K=\cup_n \holnk$ where each $\holnk$ is a twisted
form of a closed disc with respect to a finite Galois extension
(depending on $n$). The Galois cocycle giving the descent datum
comes out of a Lubin-Tate group $\G$ for $\ol$ and the logarithm
$\log_\G$ identifies $\holnk$ with a finite \'etale covering of a
closed disc whose degree equals $q^{en}.$ Here, $e$ and $q$ equal
the ramification index of $L/\Q_p$ and the cardinality of the
residue field of $L$ respectively. Passing to the limit in $n$
proves $\hol_K$ to be an \'etale covering of the affine line over
$K$. The results of sect. \ref{section-forms} may be applied to
$\holnk$ and yield the finiteness of $Pic(\holnk)$. Using a
spectral sequence argument combined with a vanishing result of L.
Gruson \cite{GrusonFi} we find $Pic(\hol_K)=\limi Pic(\holnk)$.
Building on ideas of A. de Jong \cite{deJongCrystalline} we show
that the zero section of the group $\hol_K$ supports a divisor
class of infinite order in $Pic(\hol_K)$. Finally, since the ring
of global sections $\Os(\hol_K)$ is a Pr\"uferian domain, Serre's
theorem from algebraic $K$-theory \cite{Bass} implies the result
on $K_0(\hol_K)$.

\vskip8pt

As explained above the points of $\hol$ parametrize the locally
analytic characters of $\ol$. Generalizing the construction of
$\hol$ M. Emerton has introduced such a character variety for any
abelian locally $L$-analytic group which is topologically finitely
generated \cite{EmertonA}. We conclude this work by briefly
explaining how the problem of determining the Picard group of
general character varieties can essentially be reduced to the case
of (copies of) $\hol$.

~\\{\it Acknowledgements.} I am indebted to Peter Schneider and
Jeremy Teitelbaum for their useful advice and comments during the
preparation of this article and for generously providing me with
some helpful private notes on their own work. Parts of this work
were written during a stay of the author at the {\it Tata
Institute of Fundamental Research}, Mumbai, supported by the
Deutsche Forschungsgemeinschaft. The author is grateful for the
support of both institutions.

\section{Preliminaries on Dedekind domains}\label{section-dedekind}

We recall some divisor theory for Dedekind domains, cf.
\cite{B-CA}, VII.\S2, thereby fixing some notation. Recall that an
integral domain is a {\it Dedekind domain} if it is noetherian,
integrally closed and every nonzero prime ideal is maximal.

\vskip8pt

Let $B$ be a Dedekind domain with field of fractions $F$. The free
abelian group $D(B)$ on the set $Sp(B)$ of maximal ideals of $B$
is called the {\it divisor group} of $B$. Given $P\in Sp(B)$ the
localization $B_P$ of $B$ at $P$ is a discrete valuation ring. Let
$v_P$ be the associated valuation. We have the well-defined group
homomorphism
\[div_B: F^\times\longrightarrow D(B)~,~~~ x\mapsto\sum_{P\in Sp(B)}
v_P(x)P\] whose cokernel
$$Cl(B):=D(B)/{\rm im~} div_B$$ is called the {\it divisor class group}.

 On the other hand, let $Pic(B)$
denote the {\it Picard group} of $B$, i.e. the group of
isomorphism classes of locally free $B$-modules of rank $1$ (with
the tensor product of $B$-modules as group law). There is a
commutative diagram of abelian groups with exact rows
\begin{equation}\label{cartier}\xymatrix{
1\ar[r]& B^\times \ar[r]^{\subseteq} \ar[d]^= &  F^\times \ar[r]
\ar[d]^= &
Cart(B) \ar[r]\ar[d]^{\iota} &Pic(B)\ar[d]^{\bar{\iota}}\ar[r] &1\\
1\ar[r]&  B^\times \ar[r]^{\subseteq} &  F^\times \ar[r]^{div_B} &
D(B) \ar[r]& Cl(B) \ar[r] & 1. }\end{equation} Here, $Cart(B)$
refers to the group of invertible fractional ideals of $B$ and the
map $F^\times\rightarrow Cart(B)$ is given by $x\mapsto xB$. The
map $\iota: Cart(B)\rightarrow D(B)$ is given by
$I\mapsto\sum_{P\in Sp(B)} v_P(I)P$ where $v_P(I)$ equals the
order of the extended fractional ideal in the discretely valued
field $Quot(B_P)$. Both $\iota$ and $\bar{\iota}$ are bijections,
e.g. \cite{WeibelK}, Cor. I.3.8.1.

Now suppose that $A\subseteq B$ is a subring which is a Dedekind
domain itself such that $A\rightarrow B$ is integral or flat.
Given $P'\in Sp(A),~ P\in Sp(B)$ with $P\cap A=P'$ let
$e(P/P')\in\N$ denote the ramification index of $P$ over $P'$.
Then $j(P'):=\sum e(P/P')P$ induces a well-defined group
homomorphism \[j: D(A)\rightarrow D(B)\] where the sum runs
through all $P\in Sp(B), P\cap A=P'$. It factors into a group
homomorphism
\[\overline{\jmath}: Cl(A)\rightarrow Cl(B),\] cf. \cite{B-CA},
VII.\S1.10 Prop. 14.

Now suppose additionally that $A\subseteq B$ is a finite Galois
extension with group $G$ and such that $Cl(B)=1$. The group $G$
acts on $Sp(B)$ and on $D(B)$. The map $j$ induces an isomorphism
$D(A)\car D(B)^G$. Furthermore, $Quot(A)=F^G$. Taking
$G$-invariants in the lower horizontal row of (\ref{cartier}) and
using Hilbert 90 yields the exact sequence
\[1\longrightarrow A^\times\longrightarrow Quot(A)^\times\longrightarrow
D(A)\stackrel{\delta}{\longrightarrow}
H^1(G,B^\times)\longrightarrow 1.\]

We obtain a canonical isomorphism
$$\bar{\delta}: Cl(A)\stackrel{\cong}{\longrightarrow}H^1(G,B^\times).$$


\section{Class groups of twisted affinoid discs}\label{section-forms}

Let $L$ be a complete non-archimedean field, i.e. a field that is
complete with respect to a specified nontrivial non-archimedean
absolute value. We assume that the reader is familiar with the
classical theory of affinoid spaces over such a field \cite{BGR}.

\vskip8pt

For any $L$-affinoid algebra $B$ we denote by $\Bc$ the subring of
power-bounded elements, by $\check{B}$ the $\Bc$-ideal of
topologically nilpotent elements and by $\tilde{B}:=\Bc/\check{B}$
the reduction of $B$. Passing to the reduction is a covariant
functor from $L$-affinoid algebras to algebras over the residue
field of $L$. If $|.|$ denotes the spectral seminorm on $B$ we
have $\Bc=\{b\in B: |b|\leq 1\}$ and $\check{B}=\{b\in B:|b|<1\}$.
Moreover, if the ring $B$ is reduced, the spectral seminorm is a
norm and defines the Banach topology of $B$.

\vskip8pt

After these preliminaries let $K/L$ be a finite Galois extension
and let $B$ be the one dimensional Tate algebra over $K$, i.e.

$$B=\{ \sum_{n\geq 0}a_nz^n, a_n\in K, |a_n|\rightarrow 0
{\rm~for~} n\rightarrow\infty \}.$$

Here, we denote the unique extension of the absolute value on $L$
to $K$ also by $|.|$. Suppose $A$ is a $L$-affinoid algebra
equipped with an isomorphism
$$A\otimes_L K \car B$$
of $K$-algebras. In other words, $A$ is a {\it twisted form} of
$B$ with respect to the extension $K/L$, cf. \cite{KnusOjanguren},
II.\S8. Let $G:=Gal(K/L)$. The Galois group $G$ acts on
$A\otimes_L K$ by $\sigma(a\otimes x)=a\otimes\sigma(x)$ and, via
transport of structure, we obtain in this way a semilinear Galois
action on the $K$-algebra $B$.

\vskip8pt

Remark: The forms of the affinoid algebra $B$ with respect to the
Galois extension $K/L$ are classified by the nonabelian Galois
cohomology $H^1(G, {\rm Aut}_K(B))$ according to
\cite{SchmidtDISC}. Here, ${\rm Aut}_K(B)$ refers to the
automorphism group of the $K$-algebra $B$. Any automorphism of $B$
is completely determined by its value on the variable $z$ and
induces, by functoriality, an automorphism of the algebra
$\tilde{B}=\tilde{K}[z]$. Consequently, the map $f\mapsto f(z)$
induces a group isomorphism between ${\rm Aut}_K(B)$ and the group
of formal power series
$$a_0+a_1z+a_2z^2+...$$ subject to the conditions $|a_0|\leq 1,
|a_1|=1, |a_i|<1 {\rm~for~all~}i>1$, cf. \cite{BGR}, Corollary
5.1.4/10. By the enormous size of this group the classification of
forms of $B$ seems to be a difficult task. In \cite{SchmidtDISC}
it is shown that any form of $B$ with respect to a finite tamely
ramified extension is trivial, in the sense that it is itself
isomorphic to a closed disc. Moreover, it is shown that there are
plenty of wildly ramified forms. Apart from these results the
author does not know of any results in this direction in the
literature.

\begin{lem}\label{lem-dedekind}
The ring $A$ is a Dedekind domain.
\end{lem}
\begin{proof}
The ring $B$ is a principal ideal domain. Since the extension
$A\subseteq B$ is integral the usual Going Up theorem shows $A$ to
be of dimension $1$, e.g. \cite{Matsumura}, Ex. 9.2. Finally, let
$a\in Quot(A)$ be integral over $A$ . Since $B$ is integrally
closed we have $a\in B$ and thus, by Galois invariance, $a\in A$.
\end{proof}
We therefore have the canonical isomorphism $Cl(A)\car
H^1(G,B^\times)$ at our disposal.

\begin{lem}\label{lem-sections}
Suppose that the semilinear action of $G$ maps the principal ideal
$(z)$ of $B$ to itself. Then $\Bc^\times\subseteq B^\times$
induces a short split exact sequence
\[1\longrightarrow |K^\times|/|L^\times|\longrightarrow H^1(G,\Bc^\times)\longrightarrow
H^1(G,B^\times)\longrightarrow 1.\]
\end{lem}
\begin{proof}
Let $U:=\{ b_0+b_1z+...\in B^\times: b_0=1\}$. Then
$B^\times=UK^\times$ and $\Bc^\times=U\Kc^\times$ and the Galois
action on $B$ respects these direct products by assumption. Since
cohomology commutes with direct sums we obtain, again by Hilbert
90, a short exact sequence
\[1\longrightarrow H^1(G,\Kc^\times)\longrightarrow H^1(G,\Bc^\times)\longrightarrow
H^1(G,B^\times)\longrightarrow 1.\] The map $UK^\times\rightarrow
U\Kc^\times, (u,x)\mapsto (u,1)$ induces a splitting of this
sequence. Finally, the long exact cohomology sequence induced by
the short exact sequence
$$1\longrightarrow \Kc^\times\longrightarrow
K^\times\stackrel{|.|}{\longrightarrow} |K^\times|\longrightarrow
1$$ induces, again by Hilbert 90, the isomorphism $
|K^\times|/|L^\times|\car H^1(G,\Kc^\times)$.
\end{proof}
For any real number $0<m<1$ we let $$\Bcm:=\{b\in\Bc: |b-1|\leq
m\}.$$ It is a Galois stable subgroup of $\Bc^\times$.

\begin{lem}\label{lem-sen}
There is a real number $m=m(K)$ with $0<m<1$ such that the image
of the map
\[H^1(G, \Bcm)\longrightarrow H^1(G, \Bc^\times)\]
induced by the inclusion $\Bcm\subseteq \Bc^\times$ is the trivial
group $\{1\}$.
\end{lem}
\begin{proof}
This is a mild generalization of a lemma in \cite{ST2b}. Denote by
$Tr$ the trace map of the field extension $K/L$. Since $K/L$ is
separable there is an element $c\in K$ with $Tr(c)=1$. Then any
$0<m<1$ such that $m\,|c|<1$ will do. Indeed, let $\psi$ be a
cocycle representing an element of $H^1(G, \Bcm)$ and consider the
element
\[\phi:=\sum_{\sigma\in G}\psi(\sigma)\sigma(c)\in
B.\] We have
\[|\phi-1|=|\sum_{\sigma\in
G}(\psi(\sigma)-1)\sigma(c)|\leq m\,|c|<1\] and therefore $\phi\in
\Bc^\times$. Given $\tau\in G$ we compute
$$\tau(\phi)=\sum_{\sigma\in G}\psi(\sigma)^\tau(\tau\sigma)(c)=\sum_{\sigma\in G}\psi(\tau)^{-1}\psi(\tau\sigma)(\tau\sigma)(c)=
\psi(\tau)^{-1}\phi$$ and, thus, $\psi(\tau)=\phi^{1-\tau}$.
Hence, the image of the class of $\psi$ in $H^1(G, \Bc^\times)$
coincides with the class of the coboundary $\tau\mapsto
\phi^{1-\tau}$.
\end{proof}
In the situation of the lemma the canonical homomorphism
\begin{equation}\label{injective}H^1(G,\Bc^\times)\longrightarrow H^1(G,\Bc^\times/\Bcm)\end{equation} is
therefore {\it injective}.

\vskip8pt

We fix an algebraic closure $K\subseteq \bar{K}$ of the field $K$
and we extend the absolute value from $K$ to $\bar{K}$. Let $r\in
|\bar{K}^\times |$ and consider the generalized Tate algebra $B_r$
of all affinoid functions on the closed disc of radius $r$ around
zero. It is given by all formal series

$$a_0+a_1z+a_2z^2+...$$ subject to $a_n\in K,
|a_n|r^{n}\rightarrow 0$ for $n\rightarrow\infty$. Suppose $r'\in
|\bar{K}^\times |$ is a second radius with $r<r'$. It is
well-known that, in case the field $K$ is locally compact, the
canonical inclusion
$$f: B_{r'}\stackrel{\subset}{\longrightarrow} B_r$$
is a compact continuous linear map between $K$-Banach spaces, e.g.
\cite{NFA}, Example \S16. By loc.cit., Remark 16.3 this implies
that the image of any bounded $\Kc$-module in $B_{r'}$ has compact
closure. We also remark that, since the inclusion
$Sp(B_r)\rightarrow Sp(B_{r'})$ identifies the source with an
affinoid subdomain in the target, the ring extension underlying
$f$ is flat, cf. \cite{BGR}, Cor. 7.3.2/6.

\vskip8pt

We know place ourselves in the following situation. Suppose
$r_1<r_2$ are two numbers in $|K^\times|$, the value group of $K$.
The corresponding affinoid algebra on the closed disc with radius
$r_i$ is denoted by $B_i$. We {\it suppose} that there is a
$K$-semilinear $G$-action on $B_i$ stabilizing the principal ideal
$(z)$ such that the inclusion map $f: B_2\rightarrow B_1$ is
equivariant. We suppose furthermore, that the inclusion of the
ring of invariants $A_i:=(B_i)^G$ (an affinoid $L$-algebra) into
$B_i$ is a finite Galois extension with group $G$ and that the
induced morphism between affinoids $Sp(A_1)\rightarrow Sp(A_2)$ is
an inclusion making the source an affinoid subdomain in the
target. In this situation we prove the
\begin{prop}\label{prop-finite} If the field $K$ is locally compact, the class group $Cl(A_1)$ is a finite group.
\end{prop}
\begin{proof}
First, since $r_i$ is in the value group of the ground field $K$
it is almost obvious that the results proved above for the
ordinary Tate algebra $B$ hold verbatim for the algebras $B_i$.

To start with the ring extension $A_2\subseteq A_1$ is flat and so
we have the homomomorphism $j:D(A_2)\rightarrow D(A_1)$ from the
previous section. On the other hand, the map $Sp(A_1)\rightarrow
Sp(A_2)$ induces a group homomorphism $\varphi: D(A_1)\rightarrow
D(A_2)$. According to \cite{BGR}, Prop. 7.2.2/1 (iii) it is a
splitting of $j$, i.e. $j\circ\varphi= id|_{D(A_1)}$. We have
similar maps $j$ associated to the flat extensions $A_i\subseteq
B_i$ and $f$. By multiplicativity of the ramification index they
assemble to the commutative diagram
\[\xymatrix{
D(A_2) \ar[d] \ar[r] &  D(B_2)  \ar[d] \\
D(A_1) \ar[r]  & D(B_1). }
\]
It follows that the maps $\delta_i: D(A_i)\rightarrow
H^1(G,B^\times_i)$ sit in the commutative diagram

\[\xymatrix{
D(A_2) \ar[d]^{j} \ar[r]^<<<<{\delta_2} &   H^1(G,B^\times_2) \ar[d]^{H^1(f)} \\
D(A_1) \ar[r]^<<<<{\delta_1}  &H^1(G,B^\times_1) }
\]
and, consequently, we have the identity
\begin{equation}\label{commute}H^1(f)\circ\delta_2\circ\varphi=\delta_1\circ
j\circ\varphi=\delta_1.\end{equation}

After this preliminary discussion we consider for $0<m<1$ the
following diagram $(+)$ of abelian groups

\[\xymatrix{
D(A_1) \ar[d]^{\delta_2\circ\varphi} \ar[r]^{\delta_1} &  H^1(G,B_1^\times)  \ar[d]^{=} \\
H^1(G,B_2^\times) \ar[d]^{s_2} \ar[r]^{H^1(f)} & H^1(G,B_1^\times)   \ar[d]^{s_1}  \\
H^1(G,\Bc_2^\times) \ar[d] \ar[r]  & H^1(G,\Bc_1^\times) \ar[d] \\
H^1(G, \Bc_2^\times/(\Bc_2^\times\cap\Bcm_1))\ar[r] & H^1(G,
\Bc_1^\times/\Bcm_1). }
\]
Here, the maps $s_i$ are the canonical sections from Lemma
\ref{lem-sections} and the remaining maps are the obvious ones. We
prove in a first step that all squares in this diagram are
commutative. The identity (\ref{commute}) means that the upper
square commutes. Consider the square in the middle. Since $r_i\in
|K^\times|$ we have direct product decompositions
$B_i^\times=U_iK^\times$ and $\Bc_i^\times=U_i\Kc^\times$ as in
the proof of Lemma \ref{lem-sections} where $U_i:=\{
b_0+b_1z+...\in B_i^\times: b_0=1\}$. These decompositions are
respected by $f$. The middle square is therefore induced by the
commutative diagram of Galois modules

\[\xymatrix{
B_2^\times \ar[d] \ar[r]^{f} &    \ar[d] B_1^\times \\
 \Bc_2^\times \ar[r] & \Bc_1^\times  }
\]
where the vertical maps are given by $(u,x)\mapsto (u,1)$. It is
therefore commutative. Finally, the commutativity of the lowest
square is clear.

We now prove the assertion of the proposition. In $(+)$ the upper
horizontal arrow $\delta_1$ factores into an isomorphism
$\bar{\delta}_1: Cl(A_1)\car H^1(G,B_1^\times)$. By
(\ref{injective}) we may adjust the number $m=m(K)$ so that the
lowest vertical arrow on the right hand side becomes injective.
Since the right-hand vertical maps are now all injective it
suffices to show that the image of the lower horizontal map
$$H^1(G,\Bc_2^\times/(\Bc_2^\times\cap\Bcm_1))\longrightarrow
H^1(G, \Bc_1^\times/\Bcm_1)$$ is finite. We claim that already the
set $\Bc_2^\times/(\Bc_2^\times\cap\Bcm_1)$ is finite which yields
the claim according to \cite{SerreL}, Cor. VIII.\S2.2. It is here
where we use the local compactness of the field $K$. Indeed,
consider the composite homomorphism
$$h: \Bc_2^\times\stackrel{h_1}{\longrightarrow}\Bc_1^\times\stackrel{h_2}{\longrightarrow}\Bc_1^\times/\Bcm_1$$
with kernel $\Bc_2^\times\cap\Bcm_1$. Here, $h_1$ is induced from
the compact map $f$ and $h_2$ is the canonical projection. Since
$B_i$ is reduced the spectral norm induces its Banach topology and
so $\Bc_i$ is open and closed in $B_i$. We equip $\Bc_i^\times$
with the induced topology from $\Bc_i$. The maps $h_i$ are then
continuous. We have $\overline{h_1(\Bc_2^\times)}\subseteq
\overline{f(\Bc_2)}$ and the right hand side is a compact subset
of $\Bc_1$. Hence so is the left hand side. But the units
$\Bc_1^\times$ are a closed subset of the complete adic ring
$\Bc_1$ and therefore
$\overline{h_1(\Bc_2^\times)}\subseteq\Bc_1^\times$. The image
$h_2(\overline{h_1(\Bc_2^\times)})$ is compact and contains the
image of $h$. But $\Bcm_i$ is open in $\Bc_i^\times$ and therefore
the target of $h$ is a discrete space. Hence the image of $h$ must
be finite. This proves the assertion.
\end{proof}


\section{$p$-adic Fourier theory}\label{section-fourier}
For the basic theory of locally analytic groups over $p$-adic
fields we refer to P. Schneider's monograph
\cite{SchneiderBookLie}.\vskip8pt

Let $|.|$ be the absolute value on $\C_p$ normalized via
$|p|=p^{-1}$. Let $$\Q_p\subseteq L\subseteq\C_p$$ be a finite
extension field. Let $e=e(L/\Q_p)$ be the ramification index, $k$
be the residue field of $L$ and $q=\# k$ its cardinality. Let
$\ol\subseteq L$ be the integers in $L$ and let
$e_1,...,e_{[L:\Q_p]}$ be a $\Z_p$-basis of $\ol$. We always view
$o$ as an abelian locally $L$-analytic group of dimension one.

We denote by $\B^{s}$ the rigid $L$-analytic open unit disc around
zero of dimension $s\geq 1$ and by $\B^{s}(r)$ a closed subdisc of
a real radius $r$. If $s=1$ we usually omit it from the notation.
Let $z_1,...,z_{[L:\Q_p]}$ be a set of parameters on the disc
$\B^{[L:\Q_p]}$.

\vskip8pt Let $\log (1+Z)=Z-Z^2/2+Z^3/3-...$ be the usual
logarithm series. The central object of our investigations will be
the closed analytic subvariety
$$\hol\subseteq \B^{[L:\Q_p]}$$ defined by the equations
\[ e_i\log(1+z_j)-e_j\log(1+z_i)=0\]
for $i,j=1,...,[L:\Q_p]$. It is a connected smooth one dimensional
rigid $L$-analytic variety, cf. \cite{ST2},\cite{Teitelbaum}. As
explained in the introduction it is the central object of the
$p$-adic Fourier theory developed in the article \cite{ST2}. A
particular feature of $\hol$ is that for any intermediate complete
field $L\subseteq K\subseteq \C_p$ the $K$-valued points
$z\in\hol(K)$ are in natural bijection with the set of $K$-valued
locally analytic characters $\kappa_z$ of $o$. This makes $\hol$ a
group object.

\bigskip As explained in the introduction we propose to study the
Picard group and the Grothendieck group of
$$\hol_K:=\hol\;\hat{\otimes}_L K$$ for any
intermediate field $L\subseteq K\subseteq \C_p$ which is a finite
extension of $L$. We remark straightaway that \cite{ST2}, Lem.
3.10 implies that the ideal sheaf corresponding to the zero
section in the group $\hol_K$ is an invertible sheaf whose class
in $Pic(\hol_K)$ is nontrivial if $L\neq\Q_p$. Since $\hol_K=\B_K$
in case $L=\Q_p$ and $Pic(\B_K)=1$, cf. \cite{Lazardzeros}, Thm.
7.2, this shows
\[Pic(\hol_K)\neq 1\Longleftrightarrow L\neq\Q_p.\]

Our approach to the Picard group of $\hol_K$ rests upon the fact
that $\hol$ is a twisted form of $\B$ with respect to the
extension $L\subseteq\C_p$ and that the Galois cocycle giving the
descent datum comes out of a Lubin-Tate group for $\ol$. We give
more details in the next subsection.

\bigskip

\subsection{Lubin-Tate groups}\label{additivegroup}
For a quick introduction to Lubin-Tate theory we suggest
\cite{LangCyc}. Recall that a {\it Lubin-Tate formal group} for a
fixed prime element $\pi\in\ol$ is a certain one dimensional
commutative formal group $\G$ over $\ol$ of $p$-height $[L:\Q_p]$.
It comes equipped with a unital ring homomorpism $[.]:
\ol\rightarrow {\rm End}(\G)$. We will assume that $\G$ has the
property $[\pi]=\pi X+X^q\in\ol[[X]]$. This is no essential
restriction (loc.cit., Thm. 1.1/3.1).

Viewing $\G$ as a connected $p$-divisible group let $\G'$ be its
Cartier dual and $T(\G')$ be the corresponding Tate module
\cite{Tatediv}. The latter is a free rank one $\ol$-module
carrying an action of the absolute Galois group
\[G_L:=G(\bar{L}/L)\] of $L$ which is given by a continuous
character \[\tau': G_L\longrightarrow\ol^\times.\] Denote by
$\G_m$ the formal multiplicative group over $\ol$. There is a
canonical Galois equivariant isomorphism of $\ol$-modules
\begin{equation}\label{tatemodul}
T(\G')\cong {\rm
Hom}_{\ol_{\C_p}}(\G_{\ol_{\C_p}},\mathbb{G}_{m,\ol_{\C_p}})\end{equation}
where on the right-hand side, $G_L$ acts coefficientwise on formal
power series over $\ol_{\C_p}$ and the $\ol$-module structure
comes by functoriality from the formal $\ol$-module $\G$. Choose
once and for all an $\ol$-module generator $t'$ for $T(\G')$ and
denote by
\[F_{t'}(Z)=\omega Z+...\in\ol_{\C_p}[[Z]]Z\]
the corresponding homomorphism of formal groups.
According to \cite{deJongCrystalline}, Lem. 7.3.4, we may identify
$\G$ with the rigid $L$-analytic open unit disc $\B$ around zero
and the latter becomes an $\ol$-module object in this way. The
bijection between $\C_p$-valued points $z$ of $\B$ and
$\C_p$-valued locally analytic characters $\kappa_z$ mentioned in
the introduction to this section is then given by
\[\kappa_z(g)=1+\Ft([g].z)\]
for $g\in\ol$. This bijection comes in fact from an underlying
rigid $\C_p$-analytic isomorphism
\begin{equation}\label{keyiso}\kappa:
\B_{\C_p}\car\hol_{\C_p}\end{equation} and the corresponding
Galois cocycle is given by
\begin{equation}\label{cocycle}\sigma\mapsto
[\tau'(\sigma)^{-1}]\end{equation} for $\sigma\in G_L$. Here,
$[\tau'(\sigma)^{-1}]$ is viewed as an element of the automorphism
group of the algebra $\Os(\B_{\C_p})$ in the obvious way. We point
out the trivial but useful identity
\begin{equation}\label{explog}
[\tau'(\sigma)^{-1}]=\exp_\G(\tau'(\sigma)^{-1}\log_\G(Z))
\end{equation}
for any $\sigma\in G_L$ where $\exp_\G$ and $\log_\G$ are the
formal exponential and logarithm series of $\G$ respectively. Note
also that the linear coefficient $\omega$ of the power series
$F_{t'}$ is a period (in the sense of $p$-adic Hodge theory) for
the character $\tau'$ of the absolute Galois group $G_L$. Indeed,
(\ref{tatemodul}) implies $\sigma.\Ft=\tau'(\sigma)\Ft$ and
therefore
\begin{equation}\label{period}\omega^\sigma=\tau'(\sigma)\omega\end{equation} for all $\sigma\in G_L$.
\bigskip

Let
\[L\subseteq L_\infty\subseteq\bar{L}\]
be the algebraic field extension of $L$ obtained by adjoining the
$p^n$-torsion points of the $p$-divisible group $\G'$ to $L$ for
all $n\geq 1$. By the main result of \cite{ST2}, Appendix, the
period $\omega\in\C_p$ lies in the closure of $L_\infty$. By
construction the isomorphism $\kappa$ therefore descends from
$\C_p$ to this closure. If $L\neq\Q_p$ then, according to
\cite{ST2}, Lemma 3.9, the twisted form $\hol$ of $\B$ cannot be
trivialized over a discretely valued complete subfield of this
closure. We also remark that (\ref{tatemodul}) implies that
$L_\infty$ coincides with the field extension of $L$ obtained by
adjoining all torsion points of $\G$ and all $p$-power roots of
unity to $L$. Consequently, it contains wild ramification. It is
interesting in this situation to recall from our introduction that
A. Ducros has recently shown in \cite{Ducros} that any form of the
open unit disc $\B$ with respect to a tamely ramified finite field
extension is trivial in the sense that it is itself isomorphic to
an open disc. He also showed that there are plenty of nontrivial
wildly ramified forms.

\bigskip
After this review we begin our investigations by showing in a
first step that, locally, $\hol_K$ is a twisted form of a rigid
analytic group on a closed disc which admits trivializations over
{\it finite} extensions of $K$ inside $KL_\infty$. To do this,
define for $n\geq 0$ the increasing sequence of radii
\[r_n:=r^{1/q^{en}}\] where $r:=p^{-q/e(q-1)}$ and consider
\begin{equation}\label{quasistein1}\holn:=\hol\cap\B^{[L:\Q_p]}(r_n).\end{equation} If $n$ varies these
affinoids form a countable increasing admissible open covering of
$\hol$. On the other hand, each one dimensional disc $\B(r_n)$ is
an $\ol$-module object with respect to the induced Lubin-Tate
group structure coming from $\B(r_n)\subseteq\B$. In this
situation the isomorphism (\ref{keyiso}) induces for each $n$
isomorphisms between affinoids over $\C_p$
\begin{equation}\label{quasistein2}\B(r_n)_{\C_p}\car\hol_{n,\C_p}\end{equation} according to \cite{ST2}, Thm. 3.6.
Define a formal power series
\[h_n(Z):=\exp_\G((\omega_n/\omega)\log_\G(Z))\in\C_p[[Z]]Z.\]
Since $\omega$ lies in the closure of $L_\infty$ we may fix once
and for all $\omega_n\in L_\infty$ such that
\begin{equation}\label{condition}|\omega_n/\omega-1|<p^{-n}\end{equation} for all $n\geq 0$.

\begin{lem}\label{rigidanalytic}
For all $n\in\N$ the power series $h_n$ is a rigid analytic group
automorphism of $\B(r_n)_{\C_p}$.
\end{lem}
\begin{pr}
We give the details of a proof sketched in the unpublished note
\cite{ST2b}. Let $\G(X,Y)\in\ol[[X,Y]]$ be the formal group law
underlying $\G$. Using basic properties of $\exp_\G$ and
$\log_\G$, one computes that $h_n(Z)$ equals
\begin{equation}\label{group1}
\exp_\G(\G_a(\log_\G(Z),(\omega_n/\omega-1)\log_\G(Z)))=\G(Z,\exp_\G((\omega_n/\omega-1)\log_\G(Z)))\end{equation}
as formal power series over $\C_p$ where $\G_a$ denotes the formal
additive group. By \cite{ST2}, Lem. 3.2 we have
$[p^n].\B(r_n)=\B(r)$ whence
\[ p^n\log_\G(\B(r_n))=\log_\G([p^n].\B(r_n))=\log_\G(\B(r))=\B(r)\]
where the last identity follows from \cite{LangCyc}, Lem. \S8.6.4.
Hence, on $\B(r_n)$ we have the composite of the rigid analytic
functions
\begin{equation}\label{group2}
\B(r_n)_{\C_p}\stackrel{\log_\G}{\longrightarrow}p^{-n}\B(r)_{\C_p}
\stackrel{(\omega_n/\omega)-1}{\longrightarrow}\B(r)_{\C_p}
\stackrel{\exp_\G}{\longrightarrow}\B(r)_{\C_p}.\end{equation}
Using that the group law $\G$ is defined over $\ol$ it follows
that $h_n$ is a rigid analytic function on $\B(r_n)_{\C_p}$.
Applying the same reasoning to the formal inverse
$h_n^{-1}(Z)=\exp_\G((\omega/\omega_n)\log_\G(Z))$ shows that
$h_n$ is a rigid automorphism of $\B(r_n)_{\C_p}$. It is clear
from the definition of $h_n$ that it respects the Lubin-Tate group
structure on $\B(r_n)_{\C_p}$.\end{pr}

We fix once and for all a chain of finite Galois extensions
\[L:=L_0\subseteq L_1\subseteq...\subseteq L_\infty\]
of $L$ with the property $\omega_n\in L_n$.


\begin{lem}\label{descent}
The group isomorphism
\[\kappa\circ h_n:\B(r_n)\hato_L\C_p\stackrel{\cong}{\longrightarrow}\holn\hato_L\C_p\]
is already defined over the finite extension $L_n$ of $L$.
\end{lem}
\begin{pr}
Denote by $B_n$ and $A_n$ the affinoid algebras of $\B(r_n)$ and
$\holn$ respectively. We have obvious actions of the Galois group
\[G_n:=Gal(\bar{L}/L_n)\] on $\B(r_n)_{\C_p}, B_n\hato_L\C_p$ and
${\rm Aut}_{\C_p}(B_n\hato_L\C_p)$. Here, the latter refers to the
group of $\C_p$-algebra automorphisms of $B_n\hato\C_p$. Given
$\sigma\in G_n$ and a $\C_p$-valued point $z$ of $\B(r_n)$ we find
\[
h_n^{-1}(\sigma(h_n(\sigma^{-1}(z))))=\exp_\G((\omega/\omega^\sigma)\log_\G(z))=\exp_\G(\tau'(\sigma)^{-1}\log_\G(z))
=[\tau'(\sigma)^{-1}].z\] where the middle identity and the final
identity come from (\ref{period}) and (\ref{explog}) respectively.
Denoting by $h_n^\sharp$ the algebra automorphism associated to
$h_n$ it follows for each $\sigma\in G_n$ that
\[\sigma. h_n^\sharp=h_n^\sharp\circ[\tau'(\sigma)^{-1}]\] in ${\rm Aut}_{\C_p}(B_n\hato_L\C_p)$.
By (\ref{cocycle}) the cocycle \[G_L\rightarrow {\rm
Aut}_{\C_p}(B_n\hato_L\C_p),~\sigma\mapsto[\tau'(\sigma)^{-1}]\]
gives the descent datum for the twisted form $\kappa^\sharp:
A_n\hato_L\C_p\stackrel{\cong}{\longrightarrow}B_n\hato_L\C_p$.
According to the usual formalism of Galois descent (cf.
\cite{KnusOjanguren}, \S9) we may therefore conclude that the
$\C_p$-algebra automorphism
\[h_n^\sharp\circ\kappa^\sharp:A_n\hato_L\C_p\car
B_n\hato_L\C_p\] is $G_n$-equivariant. By the existence of
topological $L$-bases for $A_n$ and $B_n$, in the sense of
\cite{NFA}, Prop. 10.1, taking $G_n$-invariants and applying
Tate's theorem $\C_p^{G_n}=L_n$ (cf. \cite{Tatediv}, Prop. 3.1.8)
yields the claim.
\end{pr}

Remark: Using that $\Ft(Z)=\exp(\omega\log_\G(Z))$ as power series
over $\C_p$ (cf. \cite{ST2}, \S4) the isomorphism $\kappa\circ
h_n$ of the preceding proposition is given on $\C_p$-points $z$ of
$\B(r_n)$ via
\[ (\kappa\circ
h_n)_z(g)=\exp(g\omega_n\log_\G(z)),~g\in\ol.\]

\begin{prop}\label{prop-cocycle2}
The Galois cocycle of the twisted form $\hol_n$ of $\B(r_n)$ with
respect to the finite extension $L_n/L$ is given by
$$\sigma\mapsto \exp_\G((\omega_n/\omega_n^{\sigma})\log_\G(Z))$$
where $\sigma\in Gal(L_n/L)$.
\end{prop}
\begin{proof}

Let $\sigma\in Gal(L_n/L)$ and let $\tilde{\sigma}\in G_L$ be any
extension to $\bar{L}$. The value of the cocycle on $\sigma$
equals the element

$$(h_n^\sharp\kappa^\sharp)\sigma (h_n^\sharp\kappa^\sharp)^{-1}
\sigma^{-1}=h_n^\sharp\kappa^\sharp\sigma(\kappa^\sharp)^{-1}(h_n^\sharp)^{-1}\sigma^{-1}
=h_n^\sharp(\kappa^\sharp\tilde{\sigma}(\kappa^\sharp)^{-1}\tilde{\sigma}^{-1})
\tilde{\sigma}(h_n^\sharp)^{-1}\sigma^{-1}$$ in the automorphism
group of the $L_n$-affinoid algebra $\Os(\B(r_n))\otimes_L L_n$,
cf. \cite{KnusOjanguren},\S9. As usual, we identify this
automorphism with its value on the parameter and consequently,
with the power series

\[
\begin{array}{cl}
h_n^\sharp(\kappa^\sharp\tilde{\sigma}(\kappa^\sharp)^{-1}\tilde{\sigma}^{-1})
\tilde{\sigma}(h_n^\sharp)^{-1}(Z)    & \stackrel{(\ref{cocycle})}{=}h_n^\sharp[\tau'(\sigma)^{-1}]\tilde{\sigma}(h_n^\sharp)^{-1}(Z) \\
     &  \\
     &\stackrel{(\ref{explog}),(\ref{period})}{=}h_n^\sharp\exp_\G((\omega/\omega^{\tilde{\sigma}})\log_\G)\tilde{\sigma}(h_n^\sharp)^{-1}(Z) \\
     & \\

     &\stackrel{(*)}{=}\exp_\G((\omega_n/\omega)(\omega/\omega^{\tilde{\sigma}})(\omega^{\tilde{\sigma}}/\omega_n^{\sigma})\log_\G(Z))\\
& \\

&=\exp_\G((\omega_n/\omega_n^{\sigma})\log_\G(Z)).\\

\end{array}
\]
Here, we used the identity
$$\tilde{\sigma}.(h_n^\sharp)^{-1}(Z)=\exp_\G((\omega^{\tilde{\sigma}}/\omega_n^\sigma)\log_\G(Z))\in
Z\C_p[[Z]]$$ in $(*)$.
\end{proof}

\bigskip
We now consider for each fixed $n$ the base extension $\holnk$ of
$\holn$ to $K$. Again, if $n$ varies these affinoids give a
countable increasing open admissible covering of $\hol_K$. As a
result of our discussion each $\holnk$ is a twisted form of the
Lubin-Tate group on $\B(r_n)_K$ trivialized over $K_n:=KL_n$ by
$\kappa\circ h_n$ and the Galois cocycle giving the descent datum
is given by the preceding proposition.

\subsection{An \'etale covering}
We briefly like to indicate an alternative intrinsic
characterization of the twisted form $\hol$ and its affinoid
subdomains $\holn$. This builds on the analytic mapping properties
of the logarithm
$$\lambda:=\log_\G$$ associated to the Lubin-Tate group $\G$. It is
best formulated in the language Berkovich analytic spaces
\cite{BerkovichBook},\cite{BerkovichEtale}. It also involves the
beginnings of A.J. de Jong's theory of \'etale covering maps for
Berkovich spaces \cite{deJongET}.

Let $\bbB$ and $\bbA^1$ be the Berkovich analytic spaces over $L$
equal to the one dimensional open unit disc around zero and the
affine line respectively. The logarithm
$$\lambda: \bbB\longrightarrow\bbA^1$$
is an \'etale and surjective morphism, cf. \cite{RameroET}, Lem.
6.1.1. Let $\bar{L}_\infty$ be the closure of
$L_\infty\subseteq\C_p$ and $G_\infty:=Gal(L_\infty/L)$. We endow
the space $\bbB_{\bar{L}_\infty}$ with the semilinear Galois
action associated to the cocycle (\ref{cocycle}). Similarly, we
endow the space $\bbA^1_{\bar{L}_\infty}$ with the semilinear
Galois action associated to the cocycle
$$\sigma\mapsto \tau'(\sigma)^{-1}Z.$$
Using (\ref{explog}) it is elementary to see that with these
definitions the map $\lambda_{\bar{L}_\infty}$ becomes
equivariant.
\begin{lem}
The Galois descent of $\bbA^1_{\bar{L}_\infty}$ is canonically
isomorphic to $\bbA^1_L$.
\end{lem}
\begin{proof}
The semilinear Galois action on $\bbA^1_{\bar{L}_\infty}$ respects
the natural increasing covering by closed discs around zero. Let
$s_n$ be a family of real numbers in $|\bar{L}^\times|$ tending
towards infinity. Let $B_n$ be the affinoid algebra of
$\bbB(s_n)_{\bar{L}_\infty}$. Let $z\in B_n$ be a parameter,
$b:=\sum_{m\geq 0}a_mz^m\in B_n$ and $\sigma\in G_\infty$. By
(\ref{period}) we have
$$\sigma.(\sum_{m\geq 0}a_mz^m)=\sum_{m\geq 0}
a_m^\sigma (\omega/\omega^\sigma)^mz^m$$ and thus $\sigma.b=b$ if
and only if $a_m/\omega^m=(a_m/\omega^m)^\sigma$ for all $m\geq
0$. Consequently, $b$ is Galois invariant if and only if
$a_m/\omega^m\in L$ for all $m\geq 0$. It follows that the subring
of Galois invariants in $B_n$ is given by

$$B_n^{G_\infty}=\{\sum_{m\geq 0}c_m(\omega z)^m, c_m\in L, |c_m| (|\omega|s_n)^m\rightarrow 0
{\rm~for~} n\rightarrow\infty\}.$$ It is canonically isomorphic to
the affinoid algebra of the closed disc of radius $|\omega|s_n$.
For varying $s_n$ these isomorphism glue to a canonical
isomorphism between the descent of $\bbA^1_{\bar{L}_\infty}$ and
$\bbA^1_L$.
\end{proof}
We point out here that
$$|\omega|=p^\nu {\rm ~~~with~}\nu=\frac{1}{p-1}-\frac{1}{e(q-1)}$$
according to \cite{ST2}, Lemma 3.4b. Let $\hol^{an}$ be the Galois
descent of $\bbB_{\bar{L}_{\infty}}$. The Galois descent of the
map $\lambda_{\bar{L}_\infty}$ is a morphism
$$\hol^{an}\longrightarrow\bbA^1_L$$
between Berkovich analytic spaces over $L$. In the following we
show that it is \'etale and surjective and therefore an \'etale
covering of the affine line over $L$.

\vskip8pt

To do this we make use of the arguments in proof of
\cite{RameroET}, Lemma 6.1.1. For $n\geq 1$ we let
$$s_n:=r|\pi|^{-en}$$
and let $\bbB(s_n)\subset\bbA^{1}$ be the closed disc of radius
$s_n$. We let $E(s_n)$ be the connected component containing zero
of $\lambda^{-1}(\bbB(s_n))$ so that the induced map
$$\lambda: E(s_n)\longrightarrow \bbB(s_n)$$
is finite \'etale and surjective. As in loc.cit. one obtains that
$E(s_n)$ is the connected component containing zero of the inverse
image of $\bbB(r)$ by $[\pi^{en}]$. Indeed, since
$r<p^{-1/e(q-1)}$, the maps $\log_\G$ and $\exp_\G$ induce
mutually inverse isomorphisms of $\bbB(r)$ (\cite{LangCyc}, \S8.6
Lemma 4) so that $\exp_\G(\pi^{en}\bbB(s_n))=\bbB(r)$. Next, we
have $[\pi](Z)=\pi Z+Z^q$ and so, by the arguments given in the
proof of \cite{ST2}, Lemma 3.2 one has
$$E(s_n)=[\pi^{en}]^{-1}(\bbB(r))=\bbB(r^{1/q^{en}})=\bbB(r_n).$$ As a result of this
discussion we have a finite \'etale surjective morphism
$$\lambda: \bbB(r_n)\longrightarrow\bbB(s_n)$$
for any $n\geq 1$. Using the compatibility with $\lambda$ and
$[\pi]$ one finds that its degree equals $q^{en}$. The same
properties hold for its base change from $L$ to the field $L_n$.
We endow the space $\bbB(r_n)_{L_n}$ with the semilinear
$Gal(L_n/L)$-action associated to the cocycle from Proposition
\ref{prop-cocycle2}. Similarly, we endow the space
$\bbB(s_n)_{L_n}$ with the semilinear $Gal(L_n/L)$-action
associated to the cocycle
$$\sigma\mapsto (\omega_n/\omega_n^\sigma)(Z).$$
This makes the map $\lambda_{L_n}: \bbB(r_n)_{L_n}\rightarrow
\bbB(s_n)_{L_n}$ equivariant. The following lemma is proved in the
same way as the preceding one.
\begin{lem}
The Galois descent of $\bbB(s_n)_{L_n}$ is canonically isomorphic
to the closed disc $\bbB(t_n)$ where $t_n:=|\omega_n|s_n$.
\end{lem}
We let $\holn^{an}$ be the Galois descent of $\bbB(r_n)_{L_n}$.
The above discussion yields a finite \'etale surjective morphism
$\holn^{an}\rightarrow\bbB(t_n)$ whose degree equals $q^{en}$.
Taking the union of these maps over all $n$ we find that the
morphism
$$\hol^{an}\longrightarrow\bbA^1_L$$
is \'etale and surjective.
\subsection{Cohomology and inverse limits}\label{cohomologylimits}
To deal with the Picard group of the rigid analytic space $\hol_K$
we need to establish a general result on sheaf cohomology and
inverse limits.

\bigskip

Recall that a rigid analytic space $X$ is called {\it quasi-Stein}
if there is a countable increasing admissible open affinoid
covering $\{X_n\}_{n\in\N}$ of $X$ such that the restriction maps
$\Os(X_{n+1})\longrightarrow\Os(X_n)$ have dense image (cf.
\cite{Kiehl}, Def. 2.3).

Let $X$ be a fixed quasi-Stein space and let $\{X_n\}$ be a
defining covering. Let $Ab(X)$ be the category of abelian sheaves
on $X$.
\begin{prop}\label{universal}
For any $\Fc\in Ab(X)$ there is a short exact sequence
\[0\longrightarrow \underleftarrow{\lim}^{(1)}_n\,\Fc(X_n)
\longrightarrow H^1(X,\Fc)\longrightarrow\underleftarrow{\lim}_n
H^1(X_n,\Fc)\longrightarrow 0\] where the right surjection is
induced by the inclusions $X_n\rightarrow X$ for each $n\in\N$ and
$H^1$ denotes sheaf cohomology.
\end{prop}
\begin{pr}
Our argument is based on the ideas in \cite{SS}, Prop. 2.4. Let
$Proj_\N(Ab)$ be the category of $\N$-projective systems over the
category of abelian groups $Ab$ and let
$F:=\underleftarrow{\lim}_n$ be the projective limit viewed as an
additive functor $Proj_\N(Ab)\rightarrow Ab$. Consider the
additive functor $G: Ab(X)\rightarrow Proj_\N(Ab),~\Fc\mapsto
(\Fc(X_n))_n.$ Both $F$ and $G$ are left exact. For any $\Fc\in
Ab(X)$ there is a short exact sequence
\[0\longrightarrow (R^1F)(G\Fc)\longrightarrow R^1(FG)(\Fc)\longrightarrow (FR^1G)(\Fc)\longrightarrow
0.\] Indeed,  an injective sheaf $\mathcal{I}\in Ab(X)$ is flasque
whence the system $G(\mathcal{I})$ is surjective and therefore
$F$-acyclic. 
The five-term exact sequence associated to the Grothendieck
spectral sequence (e.g. \cite{WeibelH} Thm. 5.8.3)
\[E_2^{pq}:=(R^pF)(R^qG)(\Fc)\Rightarrow R^{p+q}(FG)(\Fc)\]
yields the desired short exact sequence since $R^2F=0$ (cf.
[loc.cit.], Cor. 3.5.4).

After this observation let $res_n: Ab(X)\rightarrow Ab(X_n)$ be
the restriction functor. It is almost immediate that we have a
cohomological $\delta$-functor $(T^{i})_{i\geq 0}:
Ab(X)\rightarrow Proj_\N(Ab)$ where
\[T^{i}:=(H^{i}(X_n, res_n\circ(.)))_n.\]
On any $G$-ringed space the restriction of an injective sheaf to
an admissible open subset remains injective (cf.
\cite{HartshorneO}, I.\S2). Hence each $res_n$ preserves
injectives. Since $Ab(X)$ has enough injectives it follows that
each $T^{i},~i>0$ is effaceable and therefore $(T^{i})_{i\geq 0}$
is universal.
It now follows that the right-derived functors of the functor $G$
are given by $(T^{i})_{i\geq 0}$ and it remains to observe that
the sheaf property yields $F\circ G=\Gamma$, the global section
functor.
\end{pr}
\subsection{Picard groups}
In this subsection we prove the structure result on the Picard
group $Pic(\hol_K)$.
\begin{prop}\label{shortexact2}
There is an exact sequence\footnote{We turn back to our earlier
convention and denote the unit element in an abelian group by
$1$.}
\[1\longrightarrow \underleftarrow{\lim}^{(1)}_n \Os(X_n)^\times
\longrightarrow Pic (X)\longrightarrow\underleftarrow{\lim}_n
Pic(X_n)\longrightarrow 1.\]
\end{prop}
\begin{pr}
Let $\Os_X$ be the structure sheaf of $X$. We apply the last
result of the preceding subsection to $\Fc=\Os_X^\times$ and use
$H^1(X,\Os_X^\times)=Pic(X)$ (cf. \cite{FresnelVanderPut}, Prop.
4.7.2).
\end{pr}

Remark: The vanishing of the $\lim^{(1)}$-term in the above
sequence depends heavily on the base field of $X$. To give an
example let for a moment $\Q_p\subseteq K\subseteq \C_p$ by an
arbitrary complete field. The open unit disc $\B_K$ over $K$ is
quasi-Stein. A defining covering is provided by the closed
subdiscs $\B(r)_K, r<1$. In this situation
$\varprojlim_r^{(1)}\Os(\B(r)_K)^\times=1$ is equivalent to $K$
being spherically complete. This is a consequence of work of M.
Lazard, cf. \cite{Lazardzeros}, Thm. 2 and Prop. 6.

\vskip8pt

Recall that the {\it dimension} of $X$ at a point $x\in X$ equals
the dimension of the local ring at $x$.
\begin{lem}\label{equivalence} Let $X$ be equidimensional
with finitely many connected components. Taking global sections
furnishes an equivalence of exact categories between vector
bundles on $X$ and finitely generated projective $\Os(X)$-modules.
In particular, $K_0(X)\cong K_0(\Os(X))$ and $Pic(X)\cong
Pic(\Os(X))$ canonically.
\end{lem}
\begin{pr}
Following \cite{GrusonFi}, Remarque V.1. one may apply {\it
mutatis mutandis} the arguments appearing in the proof of
[loc.cit.], Thm. V.1 on each connected component of $X$.
\end{pr}

\bigskip

The space $\hol_K$ is a quasi-Stein space with respect to the
covering given by the affinoids $\holnk$. This follows from the
definition (\ref{quasistein1}) and the fact that $\B^{[L:\Q_p]}$
is quasi-Stein with respect to the covering given by the
$\B^{[L:\Q_p]}(r_n)$ (cf. remark above).
\begin{prop}\label{inverselimits}
The canonical homomorphism
$$Pic(\hol_K)\car\limi Pic(\holnk)$$
is bijective.
\end{prop}
\begin{proof}
We let $A=\Os(\hol_K)$ and $A_n:=\Os(\holnk)$. By the preceding
discussion it suffices to show that the natural homomorphism
\[Pic(A)\longrightarrow \limi Pic(A_n)\]
is injective. To show this we use an adaption of the argument
given in \cite{GrusonFi}, Prop. V.3.2. Let therefore $P$ be a
projective rank $1$ module over $A$ whose class is in the kernel
of this homomorphism. Let $P_n:=P\otimes_A A_n$. The natural
restriction map $A_{n+1}\rightarrow A_n$ is injective since this
is true for its base change to $\C_p$ (cf. (\ref{quasistein2})).
Applying $P\otimes_A (\cdot)$ we have an injective map $
P_{n+1}\rightarrow P_n$ which we view as an inclusion. The first
step now is to exhibit a well-chosen generator for each free
$A_n$-module $P_n$. This is done by induction. Let $x_1$ be any
generator for $P_1$. This starts the induction. To make the
induction step suppose $x_n$ is a well-chosen generator of the
free $A_n$-module $P_n$. Let $y_{n+1}$ be an arbitrary generator
for $P_{n+1}$. Since the image of $y_{n+1}$ under the map
$P_{n+1}\rightarrow P_n$ generates
$P_n=P_{n+1}\otimes_{A_{n+1}}A_n$ we may write $y_{n+1}=a_nx_n$
with $a_n\in A_n^\times$. Let $$B_n:=\Os(\B(r_n)_{K_n}),~~~
U_n:=\{ b_0+b_1z+...\in B_n^\times: b_0=1\},~~~G_n:=G(K_n/K).$$
Passing to a finite extension of $K_n$ (if necessary) we have
$r_n\in |K_n^\times|$ and so, as in the proof of Lemma
\ref{lem-sections}, a direct product decomposition
$B_n^\times=U_nK_n^\times$. Since the $G_n$-action on
$\B(r_n)_{K_n}$ preserves the origin taking invariants gives a
direct product decomposition
\begin{equation}\label{constant}
A_n^\times=W_nK^\times.\end{equation} with $W_n:=(U_n)^{G_n}$. Let
$a_n=wv$ with $w\in W_n, v\in K^\times$. We now let
$x_{n+1}:=v^{-1}y_{n+1}$ which completes the induction step. Note
that by construction $x_{n+1}\in W_nx_n$ for all $n$.

Now $U_n=1+U'_n$ with an $o_{K_n}$-module $U'_n$ and
$W_n=1+(U'_n)^{G_n}$ with the $\ol_K$-module $(U'_n)^{G_n}$. We
see that $$V_n:=W_nx_n\subseteq P_n$$ is a convex subset of the
$L$-Banach space $P_n$. Moreover, $B_{n+1}\subseteq B_n$ implies
$U_{n+1}\subseteq U_n$ and hence $W_{n+1}\subseteq W_n$. Thus,
$V_{n+1}=W_{n+1}x_{n+1}\subseteq W_nx_n=V_n$ and so
$$V_{n+1}\subseteq V_n$$ for all $n$. Now the map
$P_{n+1}\rightarrow P_{n}$ is a compact linear map between
$L$-Banach spaces (cf. \cite{ST5}, Lem. 6.1). This implies
(\cite{NFA}, Remark 16.3) that the image of $V_{n+1}$ in $V_n$ is
relatively compact. We may therefore argue as in \cite{GrusonFi},
Prop. V.3.2 to obtain that $\cap_n V_n\neq\emptyset$. By Theorem B
for quasi-Stein spaces (\cite{Kiehl}) any nonzero element in this
intersection generates $P$. Then $P$ is free on this generator
since ${\rm Ann}_A(x)\subseteq {\rm Ann}_{A_n}(x)=0$. Hence the
class of $P$ in $Pic(A)$ is trivial.
\end{proof}

\begin{prop} The group $Pic(\holnk)=Cl(\holnk)$ is finite for all
$n\geq 1$.
\end{prop}
\begin{proof}
We pass to a finite extension of $K_{n+1}$ (if necessary) to
obtain $r_n,r_{n+1}\in |K_{n+1}^\times|$. Let $G=Gal(K_{n+1}/K)$.
We equip $\B(r_n)_{K_{n+1}}$ and $\B(r_{n+1})_{K_{n+1}}$ with the
semilinear Galois action coming from the cocycle
$$\sigma\mapsto \exp_\G((\omega_{n+1}/\omega_{n+1}^{\sigma})\log_\G(Z))$$
for $\sigma\in G$ (Lemma \ref{prop-cocycle2}). By the Lemmas
\ref{rigidanalytic} and \ref{descent} this is an action by group
automorphisms with respect to the Lubin-Tate group structure on
these discs. In particular, the origin of these discs remains
fixed under this action. We let $B_i:=\Os(\B(r_i)_{K_{n+1}})$ and
$A_i:=\Os(\hol_{K,i})$ for $i=1,2$. The canonical inclusion
$B_2\rightarrow B_1$ is thus equivariant. Since $\omega_{n+1}$
also satisfies the defining condition (\ref{condition}) for
$\omega_n$ we have $A_n=(B_n)^G$ (and of course
$A_{n+1}=(B_{n+1})^G$). Thus, taking $G$-invariants of
$B_2\rightarrow B_1$ yields the restriction map $A_2\rightarrow
A_1$. Moreover, by the very definition of $\hol_{K,n}$ the
inclusion $\hol_{K,n}\subseteq\hol_{K,n+1}$ identifies the source
with an affinoid subdomain in the target. Since $K_{n+1}$ is
locally compact we may therefore apply Prop. \ref{prop-finite}
which yields the assertion.
\end{proof}

We now prove that, in case $L\neq\Q_p$, the group $Pic(\hol_K)$ is
not a finite group. We need two auxiliary lemmas.
\begin{lem}\label{lem-codim}
Let $R$ and $R'$ be two normal $\ol_K$-algebras of topologically
finite type which are flat over $\ol_K$. Suppose $R$ is an
integral domain. Let $R\rightarrow R'$ be a finite ring extension
such that $R\otimes_{\ol_K}K\rightarrow R'\otimes_{\ol_K}K$ is a
finite free extension of degree, say, $m$. Let
$f:Spec(R')\rightarrow Spec(R)$ be the associated morphism. There
is a closed subset $V\subset Spec(R)$ with ${\rm
codim}_{Spec(R)}V\geq 2$ such that the induced morphism
$$Spec(R')\setminus f^{-1}(V)\longrightarrow Spec(R)\setminus V$$
is finite flat of degree equal to $m$.
\end{lem}

\begin{proof}
We argue along the lines of \cite{deJongCrystalline}, Lemma 7.3.2.
Since $R$ is noetherian we may consider a finite presentation of
the $R$-module $R'$
$$R^{m_1}\stackrel{\alpha}{\longrightarrow}R^{m_0}\longrightarrow
R'\longrightarrow 0.$$ We may assume here that $0\leq m_0-m\leq
m_1$. We consider the $m$-th Fitting ideal $F_m(R')$ of the
$R$-module $R'$ (\cite{Eisenbud}, III.20.2), i.e. the ideal of $R$
generated by the $(m_0-m,m_0-m)$-minors of the matrix $\alpha$.
After inverting a prime element of $\ol_K$ the $R$-module $R'$ is
generated by $m$ elements. This implies $F_m(R')\neq 0$ (loc.cit.
Cor. 20.5/Prop. 20.6). On the complement of $V:=Spec(R/F_m(R'))$
the ideal $F_m(R')$ becomes invertible. According to
\cite{BoschLuetkeII}, Lemma 3.14 the $R$-module $R'$ is therefore
locally free of rank $m$ over the open set $Spec(R)\setminus V$.
Finally, using the normality of $R$ and $R'$ it follows as in the
proof of \cite{deJongCrystalline}, Lemma 7.3.2 that ${\rm
codim}_{Spec(R)} V\geq 2$.
\end{proof}

For any rigid $K$-analytic space $X$ we denote by
$\Os(X)^0\subseteq\Os(X)$ the $K$-subalgebra consisting of
holomorphic functions which are bounded by $1$.
\begin{lem}
Let $f: X\rightarrow\B_K$ be a finite morphism of rigid $K$-analytic spaces.
 Suppose the induced homomorphism $\Os(\B_K)^0\rightarrow\Os(X)^0$ is an integral ring extension.
 Then any function $F\in\Os(X)^0$ has only finitely many zeroes on
 $X$.
\end{lem}
\begin{proof}
Let $z$ be a parameter on the disc $\B_K$. According to
\cite{deJongCrystalline}, Lem. 7.3.4, we have
$\Os(\B_K)^0=\ol_K[[z]]$, the ring of formal power series over
$\ol_K$ in the variable $z$. Let $H\in \Os(X)^0$ and consider an
equation
$$ H^m+b_{1}H^{m-1}+...+b_m=0$$
with $b_i\in\Os(\B_K)^0$. Since $\Os(\B_K)^0$ is an integral
domain we may assume $b_m\neq 0$. If $H(x)=0$ for some $x\in X$
then $b_m(f(x))=0$. By the Weierstrass preparation theorem for
$\ol_K[[X]]$, \cite{B-CA}, VII, \S3.8 Prop. 6, the power series
$b_m$ has at most finitely many zeroes. Since $f$ has finite
fibres according to \cite{BGR}, Cor. 9.6.3/6, we conclude that $H$
has at most finitely many zeroes on $X$.
\end{proof}

\begin{prop}
Let $L\neq\Q_p$. The ideal sheaf defining the zero section of the
rigid group $\hol_K$ is not a torsion element in $Pic(\hol_K)$.
\end{prop}
\begin{pr}
Abbreviate $A:=\Os(\hol_K),~B:=\Os(\B_K)$ and
$B_{\C_p}:=\Os(\B_{\C_p})$ and let $z\in B$ be a parameter. Assume
for a contradiction that the ideal sheaf in question is a torsion
element. If $I\subseteq A$ denotes the corresponding ideal of
global sections there is $1<m<\infty$ such that $$I^m=(f)$$ with
some $f\in A$. Since the trivialization $\kappa$ preserves the
origin we have $I\mapsto (z)$ via $Pic(A)\rightarrow
Pic(B_{\C_p})$ whence $(f)\mapsto (z^m)$. Now consider $f$ as a
rigid $K$-analytic map $\hol_K\rightarrow\mathbb{A}^1_K$ into the
affine line over $K$. The composite
\[\B_{\C_p}\stackrel{\kappa}{\longrightarrow}\hol_{\C_p}\stackrel{f_{\C_p}}{\longrightarrow}\mathbb{A}^1_{\C_p}\]
is then given by a power series $F$ generating the ideal $(z^m)$
of $B$ whence
\[F(z)=az^m(1+b_1z+b_2z^2+...)\] with $a\in\C_p^\times,~b_i\in
o_{\C_p}$ according to \cite{HopkinsGrossLT}, Prop. 18.7. We may
pass to a finite extension of $K$ (if necessary) and have an
element $x\in K$ with $|x|=|a|$. Passing to $x^{-1}f$ we see that
$F$ induces a rigid $\C_p$-analytic map
$\B_{\C_p}\rightarrow\B_{\C_p}$ being the union over
$\C_p$-affinoid maps $\B(r_n)_{\C_p}\rightarrow\B(r_n)_{\C_p}$.
Hence, $f$ is in fact a rigid $K$-analytic map \[f:
\hol_K\rightarrow\B_K\] equal to the union of $K$-affinoid maps
$f_n: \holnk\rightarrow\B(r_n)_K$. The latter are induced by
functions $f_n\in\Os(\holnk)$ generating $I^m\Os(\holnk)$, the
$m$-th power of the ideal defining the zero section
$Sp~K\rightarrow\holnk$.

Next, we show that the $\Os_{\B_K}$-module $f_*(\Os_{\hol_K})$ is
locally free. Denote by $f^\sharp$ and $f_n^\sharp$ the
corresponding ring homomorphisms on global sections. Fix $n$ and
let $$A_n:=\Os(\holnk),~~~~B_n:=\Os(\B(r_n)_{K_n}).$$ We pass to a
finite extension of $K_n$ (if necessary) and have an element
$a_n\in K_n$ such that $|a_n|=r_n$. Identifying $A_n\otimes_K
K_n\simeq B_n$ via the group isomorphism $\kappa\circ h_n$ the map
$$f_n^\sharp\otimes_K K_n: B_n\longrightarrow B_n$$ is given by a
power series in $B_n$ defining the $m$-th power of the zero
section $Sp~K_n\rightarrow \B(r_n)_{K_n}$. Hence
$f_n^\sharp\otimes_K K_n$ equals the map $(a_n^{-1}z)\mapsto
(a_n^{-1}z)^m\epsilon$ with suitable $\epsilon\in B_n^\times$.
Since $|a_n^{-1}z|=1$ this is an isometry with associated graded
map
\[\gor (f_n^\sharp\otimes_K K_n): \sigma(a_n^{-1}z)\mapsto
\sigma(a_n^{-1}z)^m\sigma(\epsilon)\] and $\sigma(\epsilon)\in
(\gor K_n)^\times$. Here, $\gor$ denotes the reduction functor
from $K_n$-affinoid algebras into graded $(\gor K_n)$-algebras
introduced by M. Temkin \cite{TemkinII}. Clearly, the homomorphism
$\gor (f_n^\sharp\otimes_K K_n)$ is finite free of rank $m$ on the
homogeneous basis elements
$\sigma(a_n^{-1}z)^0,...,\sigma(a_n^{-1}z)^{m-1}$. Hence,
$f_n^\sharp\otimes_K K_n$ is finite free of rank $m$ according to
\cite{LVO}, Lem. I.6.4. By faithfully flat descent $A_n$ is
therefore a finitely generated projective $B_n$-module of rank $m$
via $f_n^\sharp$. This shows the $\Os_{\B_K}$-module
$f_*(\Os_{\hol_K})$ to be a vector bundle of rank $m$.

By \cite{GrusonFi}, V.2 Remarque $3^o$ this vector bundle must be
trivial and so $f^\sharp$ induces a finite free ring extension
$B\rightarrow A$ of degree $m$. Let $A^0$ and $B^0$ be the
holomorphic functions on $\hol_K$ and $\B_K$ respectively that are
bounded above by $1$. We claim that the ring extension $f^\sharp:
B^0\rightarrow A^0$ is integral. To see this we argue along the
lines of \cite{deJongCrystalline}, Lem. 7.3.3. Let $H\in A^0$. It
satisfies an integral equation
$$T^m+b_1T^{m-1}+...+b_m=0$$
with $b_i\in B$. We consider the ring extension induced by
$f^\sharp_n$
$$R:=\Os(\B(r_n)_K)^0\stackrel{\subseteq}{\longrightarrow}
R':=\Os(\holnk)^0.$$ It is a finite extension by \cite{BGR}, Cor.
6.4.1/6. Since $\B_K$ and $\hol_K$ are normal, $R$ and $R'$ are
normal $\ol_K$-algebras of topologically finite type which are
flat over $\ol_K$. Clearly, $R$ is an integral domain. Applying
the Lemma \ref{lem-codim} we find a closed set $V\subset Spec(R)$
with ${\rm codim}_{Spec(R)}V\geq 2$ such that the extension
$R\rightarrow R'$ is finite flat of degree $d$ over
$Spec(R)\setminus V$. Consider $H$ as an element of $R'$ via the
natural restriction map $A\rightarrow \Os(\holnk)$. It then
satisfies an equation
$$T^m+b'_1T^{m-1}+...+b'_m=0$$
with $b_i'\in\Gamma(Spec(R)\setminus V,
\Os_{Spec(R)})=\Gamma(Spec(R),\Os_{Spec(R)})=R$. Comparing these
$b_i'$ to the $b_i$ above we see that they must be equal as
elements of $\Os(\B(r_n)_K)$. Consequently, each $b_i$ is of norm
$\leq 1$ on $\B(r_n)_K\subset\B_K$ for all $n$ which means $b_i\in
B^0$. This shows $f^\sharp: B^0\rightarrow A^0$ to be an integral
extension. By the preceding lemma we see that any element of $A^0$
has at most finitely many zeroes on $\hol_K$. But according to
(the proof of) \cite{ST2}, Lem. 3.9 the nonzero holomorphic
function on $\hol_K$ given on $\C_p$-valued points via
$$\kappa_z\mapsto\kappa_z(1)-\kappa_z(0)$$ is bounded above by $1$
and has infinitely many zeroes. So we have arrived at a
contradiction.
\end{pr}

As a result of the above discussion we have the following theorem.
\begin{theo}\label{pro-p}
The group $Pic(\hol_K)$ is a profinite group. In case $L\neq\Q_p$
the isomorphism class of the ideal sheaf defining the zero section
$Sp~K\rightarrow\hol_K$ is an element of infinite order in
$Pic(\hol_K)$.
\end{theo}

We briefly explain the relation of $Pic(\hol_K)$ to the
Grothendieck group $K_0(\hol_K)$. By the Lemma \ref{equivalence}
the latter coincides with $K_0(A)$ where $A=\Os(\hol_K)$.

\vskip8pt

For the following basic notions from algebraic $K$-theory we refer
to \cite{Bass}. Let $R$ be a commutative associative unital ring,
$Pic(R)$ its Picard group and $K_0(R)$ its Grothendieck group.
Mapping $1\mapsto [R]$ induces an injective group homomorphism
$\Z\rightarrow K_0(R)$ which factores through the kernel of the
determinant $\det: K_0(R)\rightarrow Pic(R)$.
Denoting by $H_0(R)$ the abelian group of continuous maps
$Spec(R)\rightarrow\Z$ we have the rank mapping ${\rm rk}:
K_0(R)\rightarrow H_0(R)$
and a surjection \begin{equation}\label{ss}{\rm rk}\oplus\det:
K_0(R)\rightarrow H_0(R)\oplus Pic(R).\end{equation} The kernel
$SK_0(R)$ consists of classes $[P]-[R^n]$ where $P$ has constant
rank, say, $n$ and $\wedge^n P\cong R$. If $R$ is noetherian of
dimension one or a Pr\"uferian domain (i.e. an integral domain
such that any finitely generated ideal is invertible) then
$SK_0(R)=0$. Indeed, in both cases Serre's theorem (loc.cit., Thm.
IV.2.5) yields that any finitely generated projective module is
isomorphic to a direct sum of an invertible and a free module (cf.
also \cite{GrusonFi}, V.2 Remarque $3^{o}$).

 By \cite{ST2}, final discussion in sect. 3, the algebra of global sections
$\Os(\hol_K)$ is a Pr\"uferian domain which implies the
\begin{prop}\label{rkplusdet}
The rank and determinant homomorphisms give a canonical
isomorphism of abelian groups $${\rm rk}\oplus\det:
K_0(\hol_K)\car \Z\oplus Pic(\hol_K).$$
\end{prop}

\vskip8pt

Remark: Let $\gor$ be the reduction functor for $K$-affinoids
introduced by M. Temkin \cite{TemkinII}. In \cite{SchmidtDISC} the
author develops a method to compute the Picard group of twisted
forms of the closed unit disc which are '$\gor$-smooth', i.e.
whose reduction is smooth over the 'graded field' $\gor K$. A
large class of such forms is given by the tamely ramified ones.
Since the field $L_\infty$ is generated by the torsion points of
the $p$-divisible group $\G'$ it is not tamely ramified over $L$.
Therefore, the form $\holn$ is generally not tamely ramified. Even
worse, the descent datum of the form $\holn$ involves the
logarithm series $\log_\G$ (Prop. \ref{prop-cocycle2}). To compute
the reduction of $\holn$ seems to require therefore a detailed
knowledge of the coefficients of $\log_\G$ which is not available.
We have therefore not been able to prove that $\holn$ is
$'\gor$-smooth'. In fact, we are rather sceptical about $\holn$
having this property.

\subsection{General character spaces}
As explained above the space $\hol$ parametrizes the locally
analytic characters of the additive group $\ol$. In the following
we explain briefly why the problem of determining the Picard group
of general character spaces essentially reduces to the case of
(copies of) $\hol$.

\vskip8pt

So consider an arbitrary locally $L$-analytic group $Z$ which is
abelian and topologically finitely generated. Let $d:={\rm dim}_L
Z$. The following generalization of the rigid analytic character
variety $\hol$ has been introduced by M. Emerton in
\cite{EmertonA}, (6.4). Let $Rig(K)$ be the category of rigid
analytic spaces over $K$. For each $X\in Rig(K)$ let $\hat{Z}(X)$
be the group of abstract group homomorphisms
$Z\rightarrow\Os(X)^\times$ with the property that for each
admissible open affinoid subspace $U\subseteq X$ the map
\[Z\rightarrow\Os(X)^\times\stackrel{res}{\longrightarrow}\Os(U)\]
is a $\Os(U)$-valued locally $L$-analytic function on $Z$. This
defines a contravariant functor \[\hat{Z}_K: Rig(K)\rightarrow
Ab\] which is in fact representable by a smooth rigid $K$-analytic
group $\hat{Z}_K$ on a quasi-Stein space. Our character group
$\hol_K$ corresponds to the case $Z=\ol$. The association
$Z\mapsto\ZK$ is a contravariant functor that converts direct
products into fibre products over $K.$

Examples: Let $\Gm$ and $\mu_n$ be the rigid $K$-analytic
multiplicative group and the rigid $K$-analytic group of roots of
unity of order $n\geq 1$ respectively. Mapping a locally analytic
character to its value on $1$ respectively on $1 {\rm~mod~} m$
induces group isomorphisms
$$\hat{\Z}_K\car \Gm,~~~~~(\widehat{\Z/m\Z})_K\car
\mu_m$$ (loc.cit.). \vskip8pt

To get a first impression of the space $\hat{Z}_K$ we look at
dimension and number of connected components. Let $\mu\subseteq Z$
be the torsion subgroup of $Z$. By \cite{EmertonA}, Prop. 6.4.1
the inclusion of the unique maximal compact open subgroup $Z_0$
into $Z$ induces a (noncanonical) isomorphism $Z_0\times\Z^r\cong
Z$ for some unique $r\geq 0$. Hence, there is a (noncanonical)
isomorphism of rigid groups $$\ZK\car \hat{Z_0}_K\times_K\Gm^r.$$
To proceed further we impose a mild condition on the group $Z_0$.
First of all, being abelian profinite $Z_0$ contains a unique open
pro-$p$-Sylow subgroup $Z_0(p)$. The torsion part $Z_0(p)^{\rm
tor}$ of the latter group is finite and a direct factor so that
$\mu$ is finite. Any complement $Z_0(p)^{\rm fl}$ (as
$\Z_p$-module) in $Z_0(p)$ to $Z_0(p)^{\rm tor}$ has finite index
in $Z_0(p)$ and hence is open according to \cite{DDMS}, Thm. 1.17.
It is therefore naturally endowed with a structure of abelian
locally $L$-analytic group. As such is has an open subgroup which
is isomorphic, as locally $L$-analytic group, to the standard
group $\ol^d$ (in the sense of \cite{B-L}, Thm. III. 7.3.4). We
{\bf assume} in the following that the torsion part in $Z_0(p)$
admits a complement $Z_0(p)^{\rm fl}$ which is isomorphic, as
locally $L$-analytic group to $\ol^d$. We fix such an isomorphism.

\vskip8pt Example: Let $\mathbb{T}$ be a linear algebraic torus
over $L$ and $Z=\mathbb{T}(L)$ its group of $L$-rational points.
Then $Z$ is topologically finitely generated. Indeed, this is
immediate for the split part of $Z$ and follows for the
anisotropic part of $Z$ by compactness, cf. \cite{BorelTits}, Cor.
\S9.4. Now assume that $\mathbb{T}$ is split over $L$ so that we
may identify $Z=(L^\times)^d$ and $Z_0=(\ol^\times)^d$. Suppose
the ramification index $e$ of $L/\Q_p$ satisfies $1>e/(p-1)$. A
possible choice for $Z_0(p)^{\rm fl}$ is given by $(1+\pi\ol)^d.$
Moreover, the usual logarithm series followed by multiplication
with $\pi^{-1}$ provides a locally $L$-analytic group isomorphism
$Z_0(p)^{\rm fl}\car \ol^d$, e.g. \cite{NeukirchI}, Prop. 5.5.

\vskip8pt

By assumption the inclusion $ \mu\times Z_0(p)^{\rm fl}\subseteq
Z$ induces a (noncanonical) isomorphism of locally $L$-analytic
groups
\[\mu\times\ol^d\times \Z^r\stackrel{\cong}{\longrightarrow} Z\]
with unique $r\geq 0$. Applying the functor $\hat{(\cdot)}_K$
gives a (noncanonical) isomorphism of rigid groups
\[
\ZK\car \hat{\mu}_K\times_K\hol^d_K\times_K\Gm^r.\]

The space $\hat{\mu}_K$ is a finite disjoint union of points
$Sp~K_i,~i=1,...,s$ (with finite field extensions $K_i/K$). Since
each $\hol_{K_i}^d$ is connected and $\mathbb{G}_{m,K_i}^r$ is
geometrically connected, their fibre product over $K_i$ remains
connected by \cite{DucrosExcellent}, Cor. 8.4. Hence,
$$\pi_0(\hat{Z}_K)=s\leq\#\mu {\rm~~~and~~~}{\rm dim~}\hat{Z}_K=d+r.$$
We conclude with the remark that, in this situation, the natural
projection morphism $\hat{\iota}:
\hat{Z}_K\rightarrow\hat{\mu}_K\times_K\hol^d_K$ induces an
isomorphism
$$Pic(\hat{Z}_K)\car Pic(\hat{\mu}_K\times_K\hol_K^d)=\oplus_{i=1,...,s} Pic(\hol^d_{K_i}).$$
Indeed, using an argument with cohomology and inverse limits
similar to subsect. \ref{cohomologylimits} one is reduced to show
that the natural projection $\holnk^d\times_K\G_{m,n}^r\rightarrow
\holnk^d$ induces an isomorphism on Picard groups for all $n$.
Here, $\G^r_{m,n}\subset\G^r_m$ is the annulus defined by
$|p|^n\leq z_i\leq |p|^{-n}, i=1,...,r$. A finite induction
reduces to the case $d=r=1$. This case follows then by properties
of the one dimensional affinoids $\holnk$ and the fact that
$Pic(\G_{m,n})=1$. We leave the remaining details to the
interested reader.

\bibliography{mybib}
\bibliographystyle{plain}

\end{document}